\documentclass[12pt]{article}
\usepackage{latexsym}
\usepackage{epsfig}
\usepackage{amsmath,amsthm,amssymb}
\usepackage[active]{srcltx}
\usepackage{enumerate}

\usepackage[usenames]{color}
\usepackage{graphicx}

\definecolor{brown}{cmyk}{0, 0.72, 1, 0.45}
\definecolor{grey}{gray}{0.5}
\def\red{\color{red}}

\def\grey{\color{grey}}

\definecolor{duck-shit}{cmyk}{0.3,0,1,0.4}
\definecolor{my-cyan}{cmyk}{1,0.4,0,0.2}
\definecolor{my-green}{cmyk}{.8,0,1,0}

\def\grey{\color{duck-shit}}

\def\mcyan{\color{black}} \def\grey{\color{black}}  \def\red{\color{black}}

\newcounter{rot}

\def\cN{{\cal N}}
\def\bcN{\overline{\cN}}

\def\a{\alpha}  \def\d{\delta} 
\def\e{\epsilon} \def\f{\phi} \def\F{{\Phi}}  \def\g{\gamma}
\def\G{\Gamma}  
 \def\th{\theta}    \def\l{\lambda}
\def\La{\Lambda}   \def\p{\pi}
\def\r{\rho}  \def\s{\sigma} 
\def\t{\tau} \def\om{\omega}  \def\Om{\Omega}

\def\cG{{\cal G}}

\def\cB{{\cal B}}
\def\cT{{\cal T}}

\def\cR{{\cal R}}
\def\cH{{\cal H}}
\def\cU{{\cal U}}

\newtheorem{theorem}{Theorem}
\newtheorem{lemma}[theorem]{Lemma}
\newtheorem{corollary}[theorem]{Corollary}

\newcommand{\proofstart}{{\bf Proof\hspace{2em}}}

\newcommand{\proofend}{\hspace*{\fill}\mbox{$\Box$}}

\def\cW{{\cal W}}
\def\cX{{\cal X}}
\newcommand{\ooi}{(1+o(1))}

\newcommand{\oomi}{(1-o(1))}

\newcommand{\brac}[1]{\left(#1\right)}

\newcommand{\sfrac}[2]{\frac{\scriptstyle #1}{\scriptstyle #2}}
\newcommand{\bfrac}[2]{\left(\frac{#1}{#2}\right)}
\def\half{\sfrac{1}{2}}
\def\cE{{\cal E}}
\newcommand{\rai}{\rightarrow \infty}
\newcommand{\ra}{\rightarrow}
\newcommand{\rat}{{\textstyle \ra}}

\newcommand{\set}[1]{\left\{#1\right\}}
\def\sm{\setminus}
\def\seq{\subseteq}

\def\E{\mbox{{\bf E}}}
\def\Var{\mbox{{\bf Var}}}
\def\Pr{\mbox{{\bf Pr}}}
\def\whp{{\bf whp}}
\def\Whp{{\bf Whp}}

\def\cF{{\cal F}}
\newcommand{\ignore}[1]{}

\newcommand{\cA}{{\cal A}}

\newcommand{\card}[1]{\left|#1\right|}

\newcommand{\beq}[1]{\begin{equation}\label{#1}}
\newcommand{\eeq}{\end{equation}}

\def\cN{{\cal N}}
\def\dist{\;\text{dist}}

\def\ole{o(1)}
\def\bd{{\bf d}}
\def\cH{{\cal H}}

\begin{document}

\makeatletter
\title{Component structure {\mcyan of the vacant set}  induced by a random walk on a random graph.}
\author{Colin Cooper\thanks{Department of  Informatics,
King's College, University of London, London WC2R 2LS, UK}\and Alan
Frieze\thanks{Department of Mathematical Sciences, Carnegie Mellon
University, Pittsburgh PA 15213, USA.
Supported in part by NSF grant DMS0753472.
}}

\maketitle \makeatother

\begin{abstract}
We consider random walks on several classes of  graphs and explore the likely structure of the
vacant set, i.e. the set of unvisited vertices. Let $\G(t)$ be the subgraph induced by the vacant set
of the walk at step $t$.
We show that for random graphs $G_{n,p}$ (above the connectivity
threshold) and for random
regular graphs $G_r,\,r\geq 3$, the graph $\G(t)$ undergoes
 a phase transition in the
sense of the well-known Erd\H{os}-Renyi
phase transition.
{\mcyan
Thus for $t\leq (1-\e)t^*$, there is  a unique giant component, plus
components of  size $O(\log n)$,  and for
$t\geq (1+\e)t^*$ all components are of size $O(\log n)$.
For $G_{n,p}$ and $ G_r$ we give the value of $t^*$, and the size of $\G(t)$.
 For $G_r$, we also give the degree sequence of $\G(t)$, the size of the giant component (if any) of $\G(t)$
and the number of tree components of $\G(t)$ of a given size $k=O(\log n)$.}
We also show that for random digraphs $D_{n,p}$ above the strong connectivity
threshold, there is a similar directed phase transition.
Thus for $t\leq (1-\e)t^*$, there is  a unique strongly connected giant component, plus
strongly connected components of  size $O(\log n)$,  and for
$t\geq (1+\e)t^*$ all strongly connected components are of size $O(\log n)$.
\end{abstract}

\section{Introduction}
The problem we consider can be described as follows.
We have a finite graph $G=(V,E)$, and a simple random walk $\cW_u$ on $G$, starting at $u\in V$.
What is the likely component structure induced by the unvisited
vertices of $G$ at step $t$ of the walk?

Initially all vertices $V$ of $G$ are unvisited or {\em vacant}.
We regard unvisited vertices as colored {\em red}.
When $\cW_u$ visits a vertex, the vertex is re-colored {\em blue}.
Let $\cW_u(t)$ denote the position of $\cW_u$ at step $t$.
Let $\cB(t)=\set{\cW_u(0),\cW_u(1),\ldots,\cW_u(t)}$ be the set of
blue vertices at the end of step $t$, and $\cR_u(t)=V\setminus \cB_u(t)$.
Let $\G(t)$ be the subgraph of $G$ induced by $\cR(t)$.
Initially $\G(0)$ is connected, unless $u$ is a cut-vertex. As
the walk continues, $\G(t)$ will shrink to the empty graph
once every vertex has been visited. We wish to determine, as far as
possible, the likely evolution of the component
structure as $t$ increases.
In this paper we will consider three models of random graphs,  with
vertex set $V=[n]$ where $[n]=\set{1,2,\ldots,n}$.
These are the random graph $G_{n,p}$ {in which each edge of $K_n$
is included independently with probability $p$,
{\grey
 the random digraph $D_{n,p}$ {in which each  edge of $K_n$
is included independently with probability $p$ in each direction,
}
and the random graph $G_r,r\geq 3$, sampled uniformly at random
from the set of all simple $r$-regular graphs.}

Because we consider random walks on random graphs, there are two sources
of error in our proofs; (i) that we sample a graph $G$ which does not have the properties
we need, or (ii) that the random walk $\cW$ does not behave in the way we require.
The set of  graphs  $\cG'$ with properties we require have measure $\oomi$ of the graph space $\cG$.
Some of our proofs are  for   walks $\cW$ on a fixed graph  $G$ from the subset $\cG'$.
In this case, walks $\cW'$ on $G$ with the properties we require have measure $\oomi$ of $\cW$.
In other proofs, we use the method of deferred decisions, and
reveal only the parts of the graph traversed by the walk. In that way,
 the vacant set $\cR(t)$ of the walk induces a random graph, whose properties we can analyze.

Apart from $O(\cdot),o(\cdot),\Om(\cdot)$  as a  function of $n \rai$, where $n=|V|$,
we use the following:
We say $A_n\ll B_n$ or $B_n\gg A_n$ if $A_n/B_n\to 0$ as $n\to\infty$, and $A_n \sim B_n$ if $\lim_{ n \rai} A_n/B_n=1$.
The notation $\om(n)$ describes a function tending to infinity as $n \rai$.
We use the expression {\em with high probability}, \whp,  to mean, with probability
$1-o(\log^{-K}n)$ for any positive constant $K$. The variable $n$ is the size of the vertex
set of the graph, and we measure both walk and graph probabilities in terms of this.
Usually,  it will be clear that we are discussing the graph, or the walk,
but if we wish to stress this point we write \whp${_G}$ or \whp$_{\cW}$.
In the case where we use deferred decisions, if $|\cR(t)|=N$, the \whp\ statements are asymptotic in  $N$,
and we assume $N(n) \rai$.
The statement of theorems in this section  uses the annealed probability measure (graph and walk),
i.e. \whp\ relative to both graph sampling and walks.

We recall the typical evolution of the random graph
$G_{n,p}$ as $p$ increases from 0 to 1. Initially it consists of isolated vertices.
As we increase $p$ or equivalently add random edges, we find that when $p=c/n$, the maximum
component size is
logarithmic  for $c<1$, followed by a phase transition around the critical value $c=1$.
When $c>1$ the  maximum
component size is linear in $n$, and all other components have logarithmic
size.
Our aim in this paper is to show
that \whp\ $\G(t)$ undergoes a reversal of this. Thus $\G(0)$
is connected and $\G(t)$  starts to break up as $t$ increases.
There is a critical value $t^*$ such that if $t<t^*$ by a sufficient
amount then $\G(t)$ consists of a unique giant
component plus components of size $O(\log n)$. Once we
have passed through the critical value, i.e. $t>t^*$
 by a sufficient amount,
then all components are of size $O(\log n)$. As $t$ increases, the maximum component
size  shrinks to zero.
We make the following definitions. A graph with vertex set $V_1$ is {\em sub-critical} if its maximum
component size is $O(\log n)$, and
{\em super-critical} if is has a {\em unique}
component $C_1(t)$ of size $\Omega(|V_1(t)|)$, where $|V_1(t)|\gg\log n$,
and
all other components are of size $O(\log n)$.

In the case of random digraphs $D_{n,p}$ the evolution is more complex. Let $q=1-(1-p)^2$.
There is a phase transition around $q=1/n$ (i.e. $p \sim 1/2n$) for the emergence of a giant
component in the underlying graph $G_{n,q}$ of $D_{n,p}$, and, around $p=1/n$ (i.e. $q \sim 2/n$)
for the emergence of a giant strongly connected component in $D_{n,p}$.
We make the following definitions. A digraph with vertex set $V_1$ is {\em directed-sub-critical} if its maximum
strongly connected component size is $O(\log n)$, and
{\em directed-super-critical} if is has a {\em unique}
strongly connected component $C_1(t)$ of size $\Omega(|V_1(t)|)$, where $|V_1(t)|\gg\log n$,
and
all other strongly connected components are of size $O(\log n)$.

{\bf Component structure of vacant set of random graphs $G_{n,p}$ and $D_{n,p}$.}
\\
We first consider $G_{n,p}$.
We assume that
$$
np={c\log n}\text{ where } (c-1)\log n\to\infty\text{ with } n, \quad\text{ and }
 c=n^{1/\omega}.
$$
{where $\om=\om(n)$ can be any function tending to infinity with $n$.}

Let
$$t_\th=n(\log\log n+(1+\th)\log c).$$
\begin{theorem}\label{th1}
Let $\e>0$ be a small constant.
Then \whp\ in $G_{n,p}$
we have (i) $\G(t)$
 is super-critical for $t\leq t_{-\e}$,
(ii) $\G(t)$ is sub-critical for $t\geq t_{\e}$.
\end{theorem}
{\mcyan
In Section \ref{Gnp}, we prove that for $t=t_\th$, \whp\ the size of
$\cR(t)$ is  $N(t_\th) = \ooi n/(c^{1+\th} \log n)$, and that $\G(t)$ behaves as the random graph
$G_{N(t_\th),p}$.
The threshold for the giant component in $G_{N,p}$ is at $Np=1$.
For $\th <0$, let $\th =-\e$. Then $N(t_{-\e})p=c^{\e}>1$, and
there is a giant component
$C_1(t_{-\e})$ of size $\Om(N(t_{-\e}))={\Omega(n/(c \log n))}$.
The bound {$c=n^{1/\om}$}, ignores larger values of $p$. On the other hand there
is not going to be a phase transition as  $p \to 1/\log n$.
}

We next consider $D_{n,p}$. Let $q=1-(1-p)^2$.
Let $\vec \G(t)$ be the directed graph induced by the vacant set of $D_{n,p}$
and let $\G(t)$ be the undirected graph induced by the vacant set of  the underlying graph $G_{n,q}$ of $D_{n,p}$.
Theorem \ref{th1}  holds unaltered for $\G(t)$, with $nq=c \log n$. For strongly connected components of $\vec \G(t)$, we have the following
theorem.
\begin{theorem}\label{th1dir}
Let $\e>0$ be a small constant. Let $np=c \log n$.
Then \whp\ in $D_{n,p}$
we have (i) $\vec \G(t)$
 is directed-super-critical for $t\leq t_{-\e}$,
(ii) $\vec\G(t)$ is directed-sub-critical for $t\geq t_{\e}$.
\end{theorem}

{\bf Component structure of vacant set of random regular graphs.}
\\
We next consider $G_{r}$ for $r\geq 3$, constant.
Let
\beq{t*}
t^*= \frac{r(r-1) \log (r-1)}{(r-2)^2}n.
\eeq
Let
\beq{Nt}
N_t=ne^{-\frac{(r-2)t}{(r-1) n}}= ne^{-\frac{t}{\r n}},
\eeq
where $\r=(r-1)/(r-2)$, and let
\beq{pt}
p_t=e^{-\frac{(r-2)t}{\r rn}}.
\eeq
\begin{theorem}\label{th2}
Let $\e>0$ be a small constant. Then \whp\ we have
\begin{enumerate}[(i)]
\item $\G(t)$
is super-critical for $t\leq (1-\e)t^*$,  and $|C_1(t)|=\Om(n)$.
Let $p=p_t$, then
$$|C_1(t)|\sim \th N_t$$
where,
\beq{ff1}
\th=1-(1-p+p\f^{1/2})^r
\eeq
and where $\f$ is the smallest positive
solution to
\beq{ff2}
\f=(1-p+p\f^{1/2})^{r-1}.
\eeq

\item $\G(t)$ is sub-critical for $t\geq (1+\e)t^*$,
\item At some time $t\sim t^*$ the maximum {component size}
in $\G(t)$ is $n^{2/3+o(1)}$.
\end{enumerate}
\end{theorem}
When $r=3$, equations \eqref{ff1}, \eqref{ff2} give  $\f=((1-p)/p)^2$ and $\th=1-(e^{t/(6n)}-1)^3$.

We can also say something about $|\cR(t)|$ and the degree sequence of $\G(t)$.

Let {$\t_0=0$ and for $k\geq 1$ let}
\beq{t0}
\t_k=n^{1-1/k},\qquad \text{ and } \qquad t_k=\frac{\r rn\log n}{{k(r-2)+r}}.
\eeq

\begin{theorem}\label{th4}
{\mcyan
Let $\e,\; \d$ be small positive constants.
Suppose that $\log\log n\ll \om=\om(n)\ll \log n$.
\begin{enumerate}[(a)]
\item \whp\ simultaneously for all $t \le (1-\e)t_0$, $|\cR(t)|=(1+o(1))N_t$. This result also holds
\whp\ for any fixed $(1-\e)t_0\leq t \le t_0-\om n$.

\item
Let $D_s(t)$ be the number of vertices of degree $s$ in $\G(t)$. Then
\whp, simultaneously for all $0\leq s\leq r$ and for all $(\t_{r-s})^{1+\d} \le
t  \le (1-\e)t_s$,
\[
D_s(t)= (1+o(1))N_t \binom{r}{s}p_t^s(1-p_t)^{r-s}.
\]
This result also holds
\whp\ for any fixed $(1-\e)t_s\leq t \le t_s-\om n$.
\end{enumerate}
}
\end{theorem}

We make the following remarks.
\begin{itemize}
{\mcyan
\item It was proved in \cite{CFReg} that the cover time
of $G_r$ is $C(G_r) \sim ((r-1)/(r-2)) n \log n =t_0$.

\item The \whp\ concentration of  e.g. $|\cR(t)|$ in Theorem \ref{th4}
holds simultaneously for all $t \le (1-\e)t_0$. The value of
$\e$ can be made arbitrarily small. We have not attempted to optimize these results.
It is also true for any given
 $t \le t_0-\om n$,
but not proven to hold for every $t$.

}

\item The range $t_{s}-\om n\leq t\leq t_{s}+\om n$ contains the
times when the number of vertices of degree $s$ is constant
in expectation and unlikely to be concentrated around its mean.
\end{itemize}

We can give some more information about the number of small components in $\G(t)$.
Again there is a gap containing the times when the expected number of such components is constant.
{\mcyan
\begin{theorem}\label{th3}
Let $\e$ be a small positive constant. Let
$$\eta(k,t)=ne^{-t/(\r n)}\;\frac{r}{k((r-2)k+2)}\binom{(r-1)k}{k-1}p_t^{k-1}(1-p_t)^{k(r-2)+2}.$$
\begin{description}
\item[(a)]
Let
$1\leq k\leq \e\log n$ and $\e n  \le t\leq (1-\e)t_{k-1}$. Let 
$N(k,t)$ denote
the number of tree components of $\G(t)$ with $k$ vertices.
Then \whp\ for any given $t$,
$N(k,t)= (1+o(1))\eta(k,t)$.
\item[(b)]
For $k$ constant, \whp, simultaneously for all $t$ in the range
$\e n  \le t\leq (1-\e)t_{k-1}$, the number of trees with $k$ vertices is concentrated
around the value $\eta(k,t)$.
This result also holds \whp\ for any given  $t$ in the range
$(\t_{k(r-2)+2})^{1+\d} \le t \le \e n,$ where $\d$ is a small positive constant.

\end{description}
\end{theorem}
Again we make some remarks.
\begin{itemize}
\item Most  small components are trees, as e.g. \whp\ $G_r$ induces at most
$O(\log n)$ cycles of constant size.

\item Theorem \ref{th3} holds  simultaneously for  $k$ constant and any $t$.
We do not claim that  Theorem \ref{th3} holds {\bf simultaneously} for all $k=O(\log n)$ and $t$.
Our proof only show this to be true for most values of $t$.

\item The intervals containing trees of size $k$ constant are nested.
Thus,  \whp\ isolated red vertices  are the first trees to appear, (at  around time
$n^{1-1/r}$), and the last to disappear, at the cover time (around   $t_0=\r n \log n$).
\end{itemize}

Finally, we state without detailed proof the results on the
vacant set arising from the use of $k$ simultaneous
random walks, $k \ge 1$ constant.
The structure of $\G_k(t)$ is equivalent to $\G(t')$ where $t'=\ooi kt$.
In particular $t^*_k=t^*/k$.
The reasons for this assertion are based on Lemma \ref{First}
and Corollary \ref{geom}. The probability that none of $k$ independent random walks
visit a vertex $v$ during steps $T,...,t$, is the $k$-th power of the probability that
a single random walk does not visit $v$ during steps $T,...,t$.

}
\subsection{Previous work}
The only previous works on this subject that we are aware of are
Benjamini and Sznitman \cite{BS},  Windisch \cite{DW}
and \v{C}erny, Teixeira and Windisch \cite{CTD}. Papers \cite{BS}, \cite{DW} deal with the component
structure of the vacant set of a random walk on a $d$-dimensional
torus. Paper  \cite{CTD} deals with random walks on $G_r$.
It shows that \whp\ $\G(t)$ is sub-critical for $t\geq (1+\e)t^*$
and that there is a unique linear size component
for $t\leq (1-\e)t^*$. They conjecture that $\G(t)$ is super-critical
for $t\leq (1-\e)t^*$ and we prove this conjecture.

Subsequent to our posting a preliminary version on the ArXiV, \v{C}erny and Teixeira \cite{CT2} have used the 
methods of this paper to give a sharper analysis of $\G(t)$ in the critical window.

\section{Uniformity}\label{uniformity}
The main idea is to realize that the graph $\G(t)$ has a simple distribution.
First consider $G_{n,p}$.
\begin{lemma}\label{lem1}\ \\
Consider a random walk on $G_{n,p}$.
Conditional on $N=|\cR(t)|$, $\G(t)$ is distributed as $G_{N,p}$.
\end{lemma}
\proofstart
This follows easily from the principle of deferred decisions. We do not have
to decide the existence or absence of edges
between vertices in $\cR(t)$ until one of them is exposed.
\proofend

Thus to prove Theorem \ref{th1} we only need high probability estimates of $|\cR(t)|$.

We next consider $G_r$. We give two structural definitions.
\begin{lemma}\label{lem2}\ \\
Consider a random walk on $G_r$. Let $\G(t)$ have vertex set $\cR(t)$ and degree sequence
$\bd=(d_{\G(t)}(v),\;v\in \cR(t))$.
Conditional on $N=|\cR(t)|$ and degree sequence $\bd$, $\G(t)$
is distributed as $G_{N,\bd}$,
the random graph with vertex set {\red$[N]$} and degree sequence $\bd$.
\end{lemma}

\proofstart
{\mcyan
Intuitively, if we condition on $\cR(t)$ and the {\em history} and if $G_1,G_2$ are graphs with
vertex set $\cR(t)$ and if they have the same degree sequence then substituting
$G_2$ for $G_1$ will not conflict with the history
i.e. every extension of $G_1$ is an extension of $G_2$ and vice-versa. We now
give a formal explanation of this.

The {\em history} of the process can be represented by a sequence $X=(X_0,\ldots,X_{t-1})$ where
$X\in [r]^t$. This sequence is to be interpreted as follows.
{\mcyan We assume the neighbours $u_1,...,u_r$
of each vertex $v$ are given in some fixed order, e.g. increasing label order, (recall that $V=[n]$).
 When at the $j$-th vertex
$v=\cW_u(j-1)$, the
walk moved to the $X_j$-th neighbour  $u_{X_j}$} of $v$, in the given order.}
The probability space for the lemma is uniform over $\cG_r\times [r]^t$,
where $\cG_r$ is the set of $r$-regular graphs on $[n]$.
Given $\om=(G,X)$, we let $\cR_\om,\cB_\om,E_\om,\bd_\om$ denote the induced values
of $\cR(t),\cB(t)$, the edges of $G$ that
are incident with $\cB(t)$ and the degree sequence of the graph $\G_\om$ induced
by $\cR(t)$. These quantities are all determined by $\om$.
For consistency, we let $d_\om(i)$ be the degree of the $i$'th vertex of $\G_\om$
in numerical order. We can in this way associate
$\G_\om$ with $G_{\bd_{\om}}$. Fix $\cR$ and a degree sequence $\bd$. Let $\cG_{\cR,\bd}$
be the set of graphs
{
with vertex set $\cR$ and
degree sequence $\bd$. For a graph $H\in \cG_{\cR,\bd}$ we let $\Omega(H)=\set{\om:\G_\om=H}$. Then
\beq{givenX}
\Pr(H\mid X)=\frac{|\Omega(H)|}{|\cG_r|r^t}.
\eeq

We argue next that if $H_1,H_2\in \cG_{\cR,\bd}$ then $|\Omega(H_1)|=|\Omega(H_2)|$.
The lemma follows from this and \eqref{givenX}.
We define a map $\f=\f_{H_1,H_2}:\Omega(H_1)\to\Omega(H_2)$ by
$\f(G,X)=((G\setminus E(H_1))\cup E(H_2),X)$ i.e. we remove the edges of $H_1$
from $G$ and replace them by the edges of $H_2$. We only have to show that
$\f(G,X)\in \Omega(H_2)$ because then $\f_{H_1,H_2}$ will
have an inverse $\f_{H_2,H_1}$. So let $\om_2=\f(\om_1)$ where $\om_1\in \Omega(H_1)$.
We have $\cB_{\om_2}=\cB_{\om_1}$ since we have not changed
$X$ and we have not changed the set of edges incident with $\cB_{\om_1}$. This also means
that $\cR_{\om_1}=\cR_{\om_2}$
and $E_{\om_1}=E_{\om_2}$ and hence $\bd_{\om_1}=\bd_{\om_2}$. By construction,
$\G_{\om_2}=H_2$ and the lemma follows.}
\proofend

Thus to prove Theorem \ref{th2} we only need high probability estimates of the
degree sequence of $\G(t)$.
The proof of Theorem \ref{th3} can in principle be derived from this, although
we do not have a simple way of doing it.
Instead we rely on a further characterization of $\G(t)$.

We use the configuration or pairing model of Bollob\'as \cite{B} and Bender and
Canfield \cite{BC}. We start with $n$ disjoint sets
of points $W_1,W_2,\ldots,W_n$ each of size $r$. We let $W=\bigcup_{i=1}^nW_i$.
A {\em configuration} $F$ is a partition
of $W$ into $rn/2$ pairs i.e. a {\em pairing}. $\Omega$ is the set of configurations.
If $F\in\Omega$ defines
an $r$-regular multi-graph $G_F=([n],E_F)$ where $E_F=\set{(i,j):\exists
\set{x,y}\in F:x\in W_i,y\in W_j}$ i.e. we contract
$W_i$ to a vertex $i$ for $i\in [n]$.

It is known that: (i) Each simple graph arises the same number of times
as $G_F$. i.e. if $G,\;G'$  are simple, {\mcyan then $|\{F:G_F=G\}|=|\{F':G_F'=G'\}|$.}
(ii) Provided  $r$ is constant, the probability
$G_F$ is simple is bounded below by a constant. Thus if $F$ is chosen
uniformly at random from $\Omega$ then any event that
occurs \whp\ for {\mcyan $F$, occurs \whp\ for $G_F$,
and hence \whp\ for $G_r$.}

Suppose now that we generate a random $F$ using a random walk on $[n]$. To do this,
we begin with
a starting vertex $i_1$, and at the start of $t$-th step
we are at some vertex $i_t$, and  have
a partition $R_t,B_t$ of $W$ into Red and Blue points respectively.
Initially, $R_1=W$ and $B_1=\emptyset$.
In addition we have a collection
$F_t$ of disjoint pairs from $W$ where $F_1=\emptyset$.

At step $t+1$ we choose a random edge incident
with $i_t$. {\mcyan   Recall that the neighbours of $i_t$ are in a fixed order indexed by $1,...,r$.
Obviously $i_t \in \cB(t)$, as it is visited by the walk,
but we treat the configuration points in $W_{i_t}$
as blue or red, depending on whether the corresponding
edge is previously traversed (blue) or not (red).
Let $x$ be chosen randomly from $W_{i_t}$. There are two cases.}

If $x\in R_t$, then the edge is unvisited, so we choose $y$ randomly from
$R_t\setminus \set{x}$.  Suppose that $y\in W_j$.
This is equivalent to moving from
$i_t\in \cB(t)$ to $i_{t+1}=j$.
{\mcyan If $j \in \cB(t)$ this is equivalent to moving
between blue vertices on a previously unvisited edge. If $j \in \cR(t)$, this is equivalent to
moving to a previously unvisited vertex.}
We update as follows.
$R_{t+1}=R_t\setminus \set{x,y}$ and $B_{t+1}=B_t\cup\set{x,y}$,
and  $F_{t+1}=F_t\cup\set{\set{x,y}}$.

If on the other hand, $x\in B_t$ then it has previously been paired with a $y\in W_j\cap B_t$
and we move from $i_t$
to $i_{t+1}=j$ without updating. We let $R_{t+1}=R_t,B_{t+1}=B_t$ and we let $F_{t+1}=F_t$.

After $t$ steps we will have constructed a random collection $F_t$ of at most
$t$ disjoint pairs from $W$.
{\mcyan $F_t$ consists of a pairing of $B_t$, and $R_t$ is unpaired.}
In principle
we can extend $F_t$ to a random configuration $F$ by adding a
random pairing of $R_t$ to it. The next lemma summarizes this discussion.
\begin{lemma}\label{walkconfig}\

\vspace{-.3in}
\begin{description}
\item[(a)] $F_t$ plus a random pairing of $R_t$ is a {uniform} random member of $\Omega$.
\item[(b)] $i\in \cR(t)$ iff $W_i\subseteq R_t$.
\end{description}
\end{lemma}

\section{Vacancy probabilities}
As in our previous papers on random walks on random graphs, we make
heavy use of Lemma \ref{First} below.
Let $P$ be the transition matrix of the walk and let
$P_{u}^{(t)}(v)=\Pr(\cW_{u}(t)=v)$ be the $t$-step transition probability. We assume the
random walk $\cW_{u}$ on $G$ is ergodic, and thus
the random walk has  stationary distribution
$\pi$, where $\pi_v=d(v)/(2m)$.

{Suppose that the eigenvalues of $P$ are $1=\l_1>\l_2\geq \cdots\geq \l_n$.
Let $\l_{\max}=\max\set{|\l_i|:i\geq 2}$. We can make $\l_2=\l_{\max}$, if necessary,
by making the chain lazy i.e. by not moving with probability 1/2 at each step. This has no
significant effect on the analysis.}
Let $\Phi_G$ be the conductance of $G$ i.e.
\beq{conduck}
\F_G=\min_{S\subseteq V,\p_S\leq 1/2}\frac{\sum_{x\in S}\p_xP(x,\bar{S})}{\p_S}.
\eeq
Then,
\begin{align}
&1-\F_G\leq \l_2\leq 1-\frac{\F_G^2}{2}\label{mix0}\\
&|P_{u}^{(t)}(x)-\pi_x| \leq (\p_x/\p_u)^{1/2}\l_{\max}^t.\label{mix}
\end{align}
A proof of this can be found for example in
Jerrum and Sinclair \cite{JS}.

\paragraph{Mixing time of $G_{n,p},\; G_r$.}
Let $T$ be  such that, for $t\geq T$
\begin{equation}\label{4}
\max_{u,x\in V}|P_{u}^{(t)}(x)-\pi_x| =O\bfrac{\min_{x \in V}\,\pi_x}{n^3}.
\end{equation}
{\mcyan
For $G=G_{n,p}$ and $np= c \log n$, $c>1$, \whp\  the conductance $\F(G)>0$ constant, and so
\beq{TGnp}
T(G_{n,p})=O(\log n).
\eeq}
For $G=G_r$, Friedman \cite{Fried} has shown that \whp\ $\l_2\leq (2\sqrt{r-1}+\e)/r\leq 29/30$, say.
In which case we can \whp\ take
\beq{Tlogn}
T(G_r)\leq 120\log n.
\eeq

If inequality \eqref{4} holds, we say the distribution of the walk is in {\em  near stationarity}.
Fix two vertices $u,v$.
 Let $h_t=\Pr(\cW_u(t)=v)$ be
the probability that  the walk  $\cW_u$ visits  $v$ at step $t$.
Let
\begin{equation}\label{Hz}
H(z)=\sum_{t=T}^\infty h_tz^t
\end{equation}
generate $h_t$ for $t\geq T$.

We next consider the returns to vertex $v$ made by a  walk $\cW_v$, starting at $v$.
Let
$r_t=\Pr(\cW_v(t)=v)$ be the probability  that the  walk returns to
$v$ at step $t = 0,1,...$. In particular note that $r_0=1$, as the walk starts on $v$.
Let
$$R(z)=\sum_{t=0}^\infty r_tz^t$$
generate $r_t$,
and
 let
\begin{equation}
\label{Qs} R_T(z)=\sum_{j=0}^{T-1} r_jz^j.
 \end{equation}
 Thus, evaluating $R_T(z)$ at $z=1$, we have $R_T(1) \ge r_0=1$.

For $t\geq T$ let $f_t=f_t(u \rat v)$ be  the probability that the
first visit made to $v$ by the walk $\cW_u$ to $v$ {\em in the period}
$[T,T+1,\ldots]$ occurs at step $t$.
 Let
$$F(z)=\sum_{t=T}^\infty f_tz^t$$
generate  $f_t$. Then we have
\begin{equation}
\label{gfw} H(z)=F(z)R(z).
\end{equation}

The following lemma gives the probability that a walk, starting from near stationarity
makes a first visit to vertex $v$ at a given step.
For proofs of the  lemma and its corollary, {see \cite{CFGiant}.}
{\mcyan The proof differs from the earlier version given in \cite{CFReg}, in that we only
consider first visits to a vertex  $v$ after the mixing time $T$.
We use the lemma to estimate  $\E |\cR_T(t)|$, the expected number of vertices unvisited after $T$,
this differs from $\E |\cR(t)|$ by at most $T$ vertices.
}

\begin{lemma}
\label{First}
Let $R_v=R_T(1)$, where $R_T(z)$ is  from (\ref{Qs}).
For some sufficiently large constant $K$, let
\begin{equation}
\label{lamby} \l=\frac{1}{KT},
\end{equation}
where $T$ satisfies \eqref{4}.
Suppose that
\begin{description}
\item[(i)]
For some constant $\th>0$, we have
$$\min_{|z|\leq 1+\l}|R_T(z)|\geq \th.$$
\item[(ii)]$T\pi_v=o(1)$ and $T\pi_v=\Omega(n^{-2})$.
\end{description}
There exists
\begin{equation}\label{pv}
p_v=\frac{\pi_v}{R_v(1+O(T\pi_v))},
\end{equation}
 such that for all $t\geq T$,
\begin{align}
f_t(u \rat v)&=(1+O(T\pi_v))\frac{p_v}{(1+p_v)^{t+1}}+O(T \pi_v e^{-\l t/2}).\label{frat}\\
&=(1+O(T\pi_v))\frac{p_v}{(1+p_v)^{t}}\qquad for\ t\geq \log^3n. \label{frat1}
\end{align}
\end{lemma}

\begin{corollary}
\label{geom}
For $t\geq T$ let $\cA_v(t)$ be the event that $\cW_u$ does not visit $v$ at steps $T,T+1,\ldots,t$.
Then, under the assumptions of Lemma \ref{First},
\begin{align}
\Pr_\cW(\cA_v(t))&=\frac{(1+O(T\p_v))}{(1+p_v)^{t}} +O(T^2\p_ve^{-\l t/2})\label{atv}\\
&=\frac{(1+O(T\pi_v))}{(1+p_v)^t}\qquad for\ t\geq \log^3n. \label{atv1}
\end{align}
We use the notation $\Pr_\cW$ when we want to emphasize that we are
dealing with the probability space of walks on $G$.
\end{corollary}

Corollary \ref{geom} gives the probability of not visiting a single vertex in time $[T,t]$. We
need to extend this  result to certain small sets of vertices.
In particular we  need to consider sets
consisting of $v$ and a subset of its neighbours $N(v)$. Let $S$ be such a subset.

Suppose now that $S$ is a subset of $V$ with $|S|=o(n)$.
By contracting $S$ to single vertex
$\g=\g(S)$, we form a graph $H=H(S)$ in which the set $S$ is replaced by
$\g$ and the edges that were contained in $S$ are contracted to loops.
The probability of no visit to $S$ in $G$ can be found (up to a multiplicative error of
$1+O(1/n^3)$) from
the probability of a first visit to $\g$ in $H$. This is the content of Lemma \ref{contract} below.

We can estimate the mixing time of a  random walk on $H$ as follows.
Note  that the conductance of $H$ is at least that of $G$.
As some subsets of vertices of $V$ have been removed by the contraction of $S$,  the set of
values that we minimise  over, to calculate the conductance of  $H$, (see \eqref{conduck}),
 is a subset of the set of values that we minimise over for $G$.
It follows that the {conductance of $H$ is bounded below by the conductance of $G$.
{\mcyan The conductance of $G$ is constant, (see discussion
below \eqref{4})}, and so using \eqref{mix0}, \eqref{mix},
we see that the } mixing time for $\cW$ in $H$ is  $O(\log n)$.

\begin{lemma}\label{contract}\cite{CFGiant}
Let $\cW_u$ be a random walk in $G$ starting at $u \not \in S$, and let $\cX_u$ be a random walk in
$H$ starting at $u \ne \g$.
Let $T$ be a mixing time satisfying \eqref{4} in both $G$ and $H$. Then
\[
\Pr(\cA_{\g}(t);H)=\Pr(\wedge_{v \in S} \cA_{v}(t);G)\brac{1 +O\bfrac{1}{n^3}},
\]
where the probabilities are those derived from the walk in the given graph.
\end{lemma}
\proofstart

Note that $m=rn/2=|E(G)|=|E(H)|$.
Let $\cW_x(j)$ (resp. $X_x(j)$)  be the position of walk $\cW_x$ (resp. $\cX_x(j)$)
 at step $j$.
Let $\G=G,H$ and let  $P_u^s(x;\G)$ be the
transition probability in $\G$, for the walk to go from $u$ to $x$ in $s$ steps.
\begin{eqnarray}
\Pr(\cA_{\g}(t);H)&=&\sum_{x\ne \g}P^{T}_{u}(x;H)\;
\Pr(X_x(s-T)\neq \g,\; T\leq s\leq t; H)
\nonumber\\
&=&\sum_{x\ne \g}\brac{\frac{d(x)}{2m}(1+O(n^{-3}))}
\Pr(X_x(s-T)\neq \g,\; T\leq s\leq t; H)\label{close}\\
&=& \sum_{x \not \in S}
\brac{P^{T}_{u}(x;G)(1+O(n^{-3}))} \Pr(\cW_x(s-T)\not \in  S,\;T\leq s\leq t;G)\label{extra1}\\
&=& \Pr(\wedge_{v \in S} \cA_{v}(t);G)(1 +O(1/n^3))\nonumber.
\end{eqnarray}
Equation \eqref{close} follows from \eqref{4}.
Equation \eqref{extra1} follows because there is a natural measure preserving
map $\f$ between walks in $G$ that start at $x\not \in S$ and avoid $S$ and walks
in $H$ that start at $x \ne \g$ and  avoid $\g$.
\proofend

\section{The evolution of $\G(t)$ in $G_{r}$}\label{Gr}
\subsection{Estimates of $R_v$}\label{def-nice}
Let
\beq{ell1}
\ell_1=\e_1\log_rn,
\eeq
for some sufficiently small $\e_1$. A cycle $C$ is {\em small} if $|C|\leq \ell_1$.
A vertex is {\em nice} if it is at distance at least $\ell_1$ from any small cycle.
Let $\cN$ denote the nice vertices and $\bcN$ denote the vertices that are not nice.

It is
straightforward to prove by first moment calculations that:
\begin{align}
&\text{\Whp\ there are at most $n^{2\e_1}$ vertices that are not nice.}\label{nice}\\
&\text{\Whp\ there are no two small cycles within distance $\ell_1$ of each other.}\label{nice1}
\end{align}
{
The results we prove are all conditional on \eqref{nice} and \eqref{nice1} holding.
This can only inflate the probabilities of unlikely events by $1+o(1)$. This includes
events defined in terms of the configuration model as claimed in Lemma \ref{walkconfig}.
For example, if a calculation shows that an event $\cE$ has probability at most $\e$ in
the configuration model, then it has probability $O(\e)$ with respect to the corresponding
subgraph of $G$ and then we only need to multiply this by $1+o(1)$ in order to
estimate the probability conditional on \eqref{nice} and \eqref{nice1}. We will continue
relying on this without further comment.
}

A vertex $v$ is {\em tree-like} to depth $k$ if $N_k(v)$ induces a tree, rooted at $v$.
Here $N_k(v)$ denotes the set of vertices at distance at most $k$ from $v$, $k\geq 1$.
Thus a nice vertex is tree-like to depth $\ell_1/2$.

\begin{lemma}\label{lemRv}\

\vspace{-.3in}
\begin{description}
\item [(a)] If $v$ is nice then
$$R_v=(1+\ole)\r\text{ where }\r=\frac{r-1}{r-2},$$
{where the $o(1)$ term is $o(\log^{-K}n)$ for any positive constant $K$.}
\item[(b)] If $v$ is not nice then
$$R_v\leq (1+\ole)\frac{r}{r-2}$$
\end{description}
\end{lemma}
\proofstart
(a) Let $H$ denote the subgraph of $G$ induced by $N_{\ell_1/2}(v)$.
This is a tree and we can embed it into an infinite $r$-regular
tree $\cT$ rooted at $v$.
Let $\cX$ be  the walk on $\cT$,
starting from $v$, and let $X_t$ be the distance of $\cX$ from
the root vertex at step $t$.

Let {$D_0=0$, and let $D_t$} be the distance from $v$ of $\cW$ in $G$ at step $t$.
We note that we can couple $\cW_v,\cX$ so that $D_t=X_t$ up until the
first time that $D_t>\ell_1/2$.

The values of $X_t$ are as follows: $X_0=0,\; X_1=1$, and
if $X_t=0$ then $X_{t+1}=1$. If $X_t >0$ then
\begin{equation}\label{pqs}
X_t = \left\{
\begin{array}{ll}
X_{t-1}- 1& \text{ with probability } q=\frac{1}{r}\\
X_{t-1}+1 & \text{ with probability } p=\frac{r-1}{r}.
\end{array}
\right.
\end{equation}
We note the following result (see e.g. \cite{Feller}), for  a
random walk on the line $=\set{0,...,a}$ with absorbing states
$\set{0,a}$, and transition probabilities $q,p$  for moves left and
right respectively. Starting at vertex $z$, the
probability of absorption at the origin 0 is
\begin{equation}\label{absorb}
\r(z,a)= \frac{ (q/p)^z-(q/p)^a}{1-(q/p)^a}\le \bfrac{q}{p}^z,
\end{equation}
provided $q \le p$.

Let $U_\infty=\set{\exists t\ge 1: X_t=0}$, i.e. the event that the particle ever
returns to the root vertex in $\cT$.
It follows from \eqref{absorb} with $z=1$ and $a=\infty$ that
\beq{inf}
\Pr(U_\infty)=\frac{1}{r-1}.
\eeq
It follows that the expected number of visits by $\cX$ to $v$ is $\frac{1}{1-\frac{1}{r-1}}=\r$.
We write
$$R_v=\sum_{t=0}^Tr_t\text{ and }\r=\sum_{t=0}^\infty\r_t$$
where $\r_t=\Pr(X_t=v)$.

Now $r_t=\r_t$ for $t\leq \ell_1/2$ and part (a) follows {once we prove that}
\beq{tails}
\sum_{t=\ell_1/2+1}^Tr_t=\ole\text{ and }\sum_{t=\ell_1/2+1}^\infty\r_t=\ole.
\eeq
The first equation of \eqref{tails} follows from
\beq{feq}
\card{r_t-\frac{1}{n}}\leq \l_{\max}^t
\eeq
where $\l_{\max}$ is the second largest eigenvalue of the walk. This follows from \eqref{mix}.

The second equation of \eqref{tails} is proved in Lemma 7 of \cite{CFReg} where it is shown that
\beq{catalan}
\sum_{t=\ell_1/2+1}^\infty\r_t\leq \sum_{2j=\ell_1/2+1}^\infty\binom{2j}{j}\frac{(r-1)^j}{r^{2j}}
\leq \sum_{2j=\ell_1/2+1}^\infty\bfrac{4(r-1)}{r^2}^j.
\eeq
{Thus
$$R_v=\r+O(T\l_{\max}^{\ell_1/2}+T/n+(8/9)^{\ell_1/2})$$}
and part (a) follows.
(b)
We next note a property of random walks on undirected graphs which follows
from results on electrical networks (see e.g. Doyle and Snell \cite{DS}).
Let $v$ be a given vertex in a graph $G$ and $S$ a set of vertices
disjoint from $v$. Let $p(G)$, {\em the escape probability},
be the probability that, starting at $v$, the walk reaches
$S$ before returning to $v$. For an unbiased random walk,
\[
p=\frac{1}{d(v) R_{EFF}},
\]
where $R_{EFF}$ is the effective resistance between $v$ and $S$ in $G$.
We assume
each edge of $G$ has resistance 1. In the notation of this paradigm,
deleting an edge corresponds to increasing the resistance of that edge to infinity.
Thus by Raleigh's Monotonicity Law, if edges are deleted from $G$
to form a sub-graph $G'$ then $R_{EFF}'\ge R_{EFF}$.
So, if we do not delete any edges incident with $v$ then
$p' \le p$.

It follows from \eqref{nice1} that $H$ becomes a tree after removing one edge.
We can remove an edge not incident
with $v$. By the above discussion on electrical resistance we see that this will
not decrease $\Pr(U_\infty^*)$,
where this is $U_\infty$ defined with respect to $\cT^*$ which is $\cT$ less one
edge, not incident with $v$. We can argue crudely that
$$\Pr(U_\infty^*)\leq \frac{1}{r}+\frac{r-1}{r}\cdot\frac{1}{r-1}=\frac{2}{r}.$$
This is because there is an $\frac{r-1}{r}$ chance of a first move to a part of
the tree that has branching
factor $r-1$ at every vertex.

Let
$U_1^*=\{\cX$ returns to the root vertex after starting at $\ell_1/2$\}.
Then, with {$F_T$} equal to the probability of a
return by $\cW_v$ to $v$ during $[1,T]$, we have
\beq{fT}
{F_T}\leq \Pr(U_\infty^*)+T\Pr(U_1^*).
\eeq
The RHS of \eqref{fT} is at least
the probability that $\cW_v$ returns before reaching distance $\ell_1/2$ or
returns after reaching distance
$\ell_1/2$ at some time $t\leq T$.

Now, using \eqref{absorb}, we see that
\beq{Hoef}
{\Pr(U_1^*)} \leq \frac{1}{(r-1)^{\ell_1/4}}.
\eeq
Here we have $\ell_1/4$ in place of $\ell_1/2$ to account for the one place
where we move left with probability
$\frac{1}{r-2}$. We argue that at least one of the paths from $v$ to $w$ or $w$
to the boundary must be
at least $\ell_1/4$ and not use the vertex incident to the deleted edge.

Thus ${F_T}\leq(2+\ole)/r$ and since $R_v\leq \frac{1}{1-F_T}$
we have $R_v\leq \frac{r+\ole}{r-2}$.
\proofend

\subsection{Proof of Theorem \ref{th4}} \label{proof of th3}

{\mcyan

This section establishes the \whp\ values of $|\cR(t)|$ and $D_s(t)$ for $\G(t)$.
We also calculate the \whp\ value of $|\cU(t)|$, where $\cU(t)$ is the
{\em number of unvisited edges at step} $t$.
This gives us the value of $|R_t|= 2|\cU(t)|$ (see Section \ref{uniformity}),
which we need for the proof
of Theorem \ref{th3}.

Fix $v\in V$ and let
$N(v)=\set{w_1,w_2,\ldots,w_r}$ and choose $0\leq s\leq r$.
Let
\[
P(v,s,t)=\Pr_\cW(\{v,w_1,\ldots,w_s\}\subseteq\cR(t) \text{ and } \{w_{s+1},\ldots,w_r\}
\subseteq \cB(t)).
\]
Then
\begin{align}
\Pr_\cW(\text{$v \in \cR(t)$ and has degree $s$ in $\G(t)$})=\binom{r}{s}P(v,s,t).\label{Pis}
\end{align}
}
Lemma \ref{First} is only valid after the mixing time $T$, and
Corollary \ref{geom}
 gives
precise results only after some $t$ sufficiently larger than $T$.
We  assume henceforth  that $t\geq  \log^3n$. We deal with
the very beginning of the walk  ($t<T$) in Section \ref{begin}.
Define the set $\cR_T(t)$ to be those vertices visited by the walk during steps $T,...,t$, and let
 $\cB_T(t)=V\setminus \cR_T(t)$. Thus
$\cR(t)=\cR_T(t) \sm  \set{\cW_u(0),\cW_u(1),\ldots,\cW_u(T-1)}$

We recall the definition of  a nice vertex, as given at the start of Section \ref{def-nice}.
We will say an
 edge $e=\{u,v\}$  is nice, if both $u,v$ are nice.
The next lemma gives  enough information to compute the  expected number of unvisited
nice vertices, edges between nice vertices, and the expected degree
sequence of nice vertices in $\cR(t)$.
\begin{lemma}\label{key} Let $p_t=e^{-((r-2)t)/(\r rn)}$ as given by \eqref{pt}.
${\bf Whp}_G$\footnote{We use the subscript $G$ to
emphasize that the probability space is random $r$-regular graphs},
for all nice vertices $v$, and for all nice edges $\{u,v\}$, for $t \ge \log^3 n$,
\vspace{-.25in}
\begin{description}
\item[(a)] $\Pr_\cW(v\in {\cR_T(t)})=(1+\ole)e^{-\frac{t}{\r n}}$.
{\mcyan \item[(b)] $\Pr_\cW(\{u,v\}\in \cU_T(t))=(1+\ole)e^{-\frac{2t}{\r rn}}$.}
\item[(c)]
Let
\[
P_T(v,s,t)=\Pr_\cW(\{v,w_1,\ldots,w_s\}\subseteq{\cR_T(t)}
\text{ and } \{w_{s+1},\ldots,w_r\}\subseteq {\cB_T(t)}),
\] then
\[
P_T(v,s,t)=
(1+\ole)e^{-\frac{t}{\r n}}p_t^s(1-p_t)^{r-s}.
\]
\end{description}
\end{lemma}
\proofstart
Part (a) follows directly from Corollary \ref{geom} using $p_v$ from \eqref{pv}
and the value of $R_v$ from \ref{lemRv}(a).

Part (b) follows similarly, by subdividing the edge $e=\{u,v\}$ with an artificial
vertex $\a$ to form the graph $H$.
The stationary distribution of $\a$ in $H$ is  $\pi_\a=2/(rn+2)$,
and the value of $R_\a$  is given by  $R_v$ of Lemma \ref{lemRv}(a).

{\mcyan
For part (c) we proceed as  follows.
Let $X \seq N(v)$.
Let $\g_X$ denote the contraction of $\set{v}\cup X$.
Corollary \ref{geom} applies to $\g_X$ and
\beq{pigx}
p_{\g_X}= \ooi  \frac{(r-2)(r+(r-2)|X|)}{r(r-1)\;n},
\eeq
as we now explain.
From \eqref{pv}, the expression for $p_{\g_X}\sim \pi_{\g_X}/R_{\g_X}$.
The stationary distribution of
 $\g_X$ is $\pi_{\g_X}=(|X|+1)/n$.
Since $v$ is nice,
the expected number of returns, is (up to a factor $1+\ole$),
$R_{\g_X}=\frac{1}{1-f}$; where $f$ is the probability of return
to the root $\g_X$ of an infinite tree with branching factor $r-1$ at
each non-root vertex. Thus $p_{\g_X}=\ooi \pi_{\g_X}(1-f)$. It remains to calculate $f$.
At the root there are
$|X|$ loops and $r-|X|+(r-1)|X|$ branching edges.
This gives
\[
f=\frac{2|X|}{r(|X|+1)}+\frac{(r-2)|X|+r}{r(|X|+1)}\frac{1}{r-1},
\]
and hence \eqref{pigx} above.}

Let $N(v)=\set{w_1,w_2,\ldots,w_r}$, and let $X=\{w_1,\ldots,w_s\} \cup Y$
where $Y \seq \set{w_{s+1},\ldots,w_r}$.
Then
\begin{align}
\Pr_\cW(&\{v,w_1,\ldots,w_s\}\subseteq{\cR_T(t)}\text{ and } \{w_{s+1},\ldots,w_r\}
\subseteq {\cB_T(t)})\nonumber\\
&=\sum_{Y\subseteq {\set{w_{s+1},\ldots,r}}}(-1)^{|Y|}\Pr_\cW((\{v,w_1,\ldots,w_s\}
\cup Y)\subseteq{\cR_T(t)})\nonumber\\
&=\sum_{Y\subseteq {\set{w_{s+1},\ldots,r}}}(-1)^{|Y|}\
\frac{1+O(T\p_{\g_X})}{(1+p_{\g_X})^t}\label{pie},
\end{align}
where $|X|=0,...,r-s$. If  $s=0$,  we suppose that $\{v,w_1,\ldots,w_s\}=\{v\}$.

Thus, for $s=0,1,...,r$,
\begin{multline}\label{nmn}
\Pr_\cW(\{v,w_1,\ldots,w_s\}\subseteq{\cR_T(t)}
\text{ and } \{w_{s+1},\ldots,w_r\}\subseteq {\cB_T(t)})=\\
\exp\set{-(1+\ole)\frac{(r-2)^2s+r(r-2)}{r(r-1)n}t}\sum_{Y\subseteq [s+1,r]}(-1)^{|Y|}
\exp\set{-(1+\ole)\frac{(r-2)^2|Y|}{r(r-1)n}t}.
\end{multline}
There are two cases.
When $t=O(n)$ we can write \eqref{nmn} as
$$\exp\set{-(1+\ole)\frac{(r-2)^2s+r(r-2)}{r(r-1)n}t}\brac{1-\exp\set{-\frac{(r-2)^2}{r(r-1)n}t}}^{r-s}
+\ole.$$
and (c) follows, since the terms above are $\Omega(1)$.
When $t/n \to\infty$ we go back to \eqref{nmn} and observe that the sum is $1-o(1)$
and thus
\begin{multline*}
\Pr_\cW(\{v,w_1,\ldots,w_s\}\subseteq{\cR_T(t)}\text{ and }
 \{w_{s+1},\ldots,w_r\}\subseteq {\cB_T(t)})=\\(1+\ole)
\exp\set{-\frac{(r-2)^2s+r(r-2)}{r(r-1)n}t}
\end{multline*}
as required.
\proofend

{\mcyan
{\bf \underline{Proof of Theorem \ref{th4}(a).}}

Let $t \le t_0-\om n$,
where $t_0=\r n \log n$ and $\om$ satisfies the conditions of Theorem \ref{th4}.

Using \eqref{nice}, \eqref{nice1} we have \whp, that $|\bcN|=O( n^{2 \e_1}\log n )$,
and thus $|\cN|=n(1-o(1))$.
Let $Z(t)=|\cR(t)\cap \cN|$. As there are at most $T=O(\log n)$ vertices
in $\cR(t)\setminus \cR_T(t)$,
\begin{eqnarray}
 \E Z(t)&=& n (1-o(1)) \Pr_\cW(v\in {\cR_T(t)})+O(T)\label{Lin1}\\
&=& (1+o(1)) n e^{-\frac{t}{\r n}} \label{Lin3}.
\end{eqnarray}
Thus $\E Z(t) \rai$
for $t \le \r n \log n - \om n$.

We  use the Chebyshev inequality to prove concentration of $Z(t)$
for a single $t \le t_0-\om n$.

Suppose that
{$\log\log n\ll \om'=\om'(n)\ll \log n$. We first show that
\beq{cheby}
\Var(Z(t))=O(r^{\om'}\E Z(t))+e^{-a\om'}\E(Z(t))^2,
\eeq
for some constant $a>0$.

Fix $t$ and let $\cE_v$ be the event that vertex $v\in \cR(t)$.
We claim that if $v,w$ are at distance at least
$\om'$ then
\beq{var}
\Pr(\cE_v\cap \cE_w)=(1+e^{-\Omega(\om')})\;\Pr(\cE_v)\Pr(\cE_w).
\eeq
To prove this we use Lemma \ref{contract}. Let $S=\set{v,w}$, and let $\g(S)$
be the contraction of $S$. For the
random walk on the associated $H$
we have
\[
\half(R_v+R_w) \le R_{\g} \le\half(R_v+R_w)(1+O(Te^{-\Om(\om')})).
\]
Indeed, the first move from $\g$ will be to a neighbour of $v$ or $w$.
Assume it is to a neighbour of $v$. The expected number of returns  will be
 $R_v$  plus $R_w$ times the probability of a visit to $w$ during the mixing time.
Because $v$ and $w$ are at distance at least $\om'$, using \eqref{feq},
 the probability of a visit to $w$ during $T$
can be bounded by $T(n^{-1}+\l_{\max}^{\om'})$.
Thus
\beq{32}
R_\g=\r\brac{1+e^{-\Omega(\om')}}, \qquad \p_\g=\frac{2}{n},
\qquad p_\g=(1-O(Te^{-\Om(\om')}))\frac{2}{\r n}.
\eeq
{Equation \eqref{var} follows on using Lemmas \ref{First}, \ref{contract} and
equation \eqref{32}.}

Thus
\begin{align*}
\E(Z^2(t))&=\E Z(t)+\sum_{v,w \atop dist(v,w)\geq \om'}\Pr(\cE_v\cap \cE_w)+
\sum_{v,w \atop dist(v,w)< \om'}\Pr(\cE_v\cap \cE_w)\\
&\leq \E Z(t)+(1+e^{-a\om'})\E(Z(t))^2+r^{\om'}\E(Z(t))
\end{align*}
and \eqref{cheby} follows.
Applying the Chebyshev inequality we see that
\beq{toterr}
\Pr \brac{|Z(t)-\E Z(t)|\geq \E(Z(t))e^{-a\om'/3}}
\leq \frac{2r^{\om'} e^{a\om'}}{\E Z(t)}+e^{-a\om'/3}.
\eeq
Provided $t \le t_0-\om n$, $\E Z(t) \ge e^{-\om/\r}/2$ and our choice of $\om'$
implies that the RHS of \eqref{toterr} is $o(1)$ for such $t$.

To prove the concentration of $Z(t)$ {simultaneously} for all $t \le (1-\e)t_0$, we proceed as follows.
Let {now} $\om'=\ell_1/2$ where $\ell_1=\e_1\log_r n$ is given by \eqref{ell1}.
From \eqref{Lin1}-\eqref{Lin3}, for $t \le (1-\e)\r n \log n$,
we have $\E Z(t) =\ooi ne^{-t/\r n}\ge n^{\e}/2$.
We see that  \eqref{toterr} becomes
\beq{toterr1}
\Pr(|Z(t)-\E Z(t)|\geq \E(Z(t))e^{-a\om'/3})
\leq \frac{2r^{\om'} e^{a\om'}e^{t/(\r n)}}{n}+e^{-a\om'/3}=o(n^{-\d}),
\eeq
for some (small) $\d>0$ constant.
We can make $\d=(a\e_1 \log r)/4$ provided we make $\e>\e_1\brac{1+\frac{a}{\log r}}$.

We interpolate $[0,t_0]$ at $M=n^{\d/2}$ integer points $s_1,...,s_M$
a distance $\s=t_0\;n^{-\d/2}$ apart (ignoring rounding). Let $\cH(M)$ be the event that
$\set{|Z(t)-\E Z(t)|\leq \E(Z(t))e^{-a\om'/3}}$
 holds simultaneously at $s_1,...,s_M$. Then $\Pr(\neg \cH(M))=o(n^{-\d/2})=o(1)$.

The value of $Z(t)$ is non-increasing, with $\E Z(t)=\ooi n e^{-t/\r n}$. Thus for any $t$,
$s_i \le t \le s_{i+1}$,
\begin{eqnarray*}
Z(s_i)& \ge &Z(t)\quad  \ge Z(s_{i+1})\\
\ooi e^{(t-s_i)/\r n} &\ge& \frac{Z(t)}{\E Z(t)} \ge \oomi e^{(t-s_{i+1})/\r n}.
\end{eqnarray*}
But $e^{(t-s_i)/\r n},e^{(t-s_{i+1})/\r n}$ are both $1-o(1)$ and so
the concentration result holds for nice vertices for $t \le (1-\e)t_0$.

Lemma \ref{key}  applies only to nice vertices.
We next consider $\bcN\cap \cR(t)$.
It follows from Lemmas \ref{First} and \ref{lemRv} that
\beq{ignore}
\whp\ \bcN\subseteq {\cB_T(t)} \text{ for } t\geq 10\e_1n\log n.
\eeq
However for $t \le 10\e_1n\log n$,
\beq{znn}
\E Z(t) \ge n^{1-10\e_1} \gg |\bcN| \ge |\bcN\cap \cR(t)|.
\eeq
This completes the proof of part (a) of Theorem \ref{th4}.

{\bf \underline{Proof of Theorem \ref{th4}(b).}}

The proof of part (b) is similar to that of (a). Observe first that
part (a) and Lemma \ref{key} imply
\beq{EDst}
\E(D_s'(t))=(1+o(1))N_t\binom{r}{s}p_t^s(1-p_t)^{r-s}.
\eeq
where $D_s'(t)$ is the number of nice vertices of degree $s$ in $\G(t)$.

Our next aim is to show
\beq{Dst}
\E(D_s'(t))=e^{\Omega(\om)}\ for\ (\t_{r-s})^{1+\d} \le t\leq t_{s}-\om n
\eeq
As a function of $t$, $N_tp_t^s(1-p_t)^{r-s}$
is log-concave and to bound $\E(D_s'(t))$ from below
we only need to check $\E(D_s'(t))$ at the lower bound $t=(\t_{r-s})^{1+\d}$ and at the upper bound
$t=t_{s}-\om n$. See \eqref{t0} for the definitions of $\t_s,t_s$.
For $t=o(n)$ we have $p_t=1-o(1)$ and $N_t\sim n$ and then \eqref{Dst} is simple to verify
at the lower bounds for $t$ in \eqref{Dst}. When $t=t_s-\om n$ we have
$p_t=e^{\frac{\om(r-2)}{\r r}}n^{-\frac{r-2}{s(r-2)+r}}$ and $N_t=e^{\frac{\om}{r}}
n^{1-\frac{r}{s(r-2)+r}}$ and \eqref{Dst} is also easy to verify.

We can use the Chebyshev inequality to prove concentration.
We let $\cF_v(s)$ be the event that $v$ is a vertex of degree $s$ in $\G(t)$. We prove that
for $v,w\in\cN$,
\beq{fkk}
\Pr(\cF_v(s)\cap \cF_w(s))=(1+e^{-\Omega(\om')})\Pr(\cF_v(s))\Pr(\cF_w(s))
\eeq
whenever $v,w$ are at least
$\om'$ apart, where
{$\om'/\log\log n\to\infty$ and $\om=o(\log n)$
and $\om'=o(\om)$}.. We can argue for this by a small change in the argument for \eqref{var}.
This proves concentration for $(\t_{r-s})^{1+\d}\leq t\leq t_{s+1}-\om n$.
Let $D_{\ge s}(t)=\sum_{k=s}^rD_k(t)$. This sum is  monotone non-increasing with $t$.
Simultaneous concentration of $D_{\geq s}(t)$ and
hence $D_s(t)$ follows from the interpolation method of part (a) applied to $D_{\ge s}(t)$.

We must now argue that the contribution of $\bcN$ is negligible.
Equation \eqref{ignore} shows that \whp\ $\bcN \seq \cB(t)$
for $t\geq T_1=10 \e_1 n \log n$.
On the other hand, by \eqref{EDst}, for $t=t^*$ and any $0 \le s \le r$,
$$
\E(D_s'(t^*))=\Omega(n^{1-\frac{10\e_1(r+s(r-2))}{\r r}})\gg n^{2\e_1}.
$$
It follows from this and concentration of $D_s'(t)$ and \eqref{nice} that $D_s'(t)\sim D_s(t)$
\whp\ and the proof of part (b) is complete.}
\ignore{
{\bf \underline{Proof of Theorem \ref{th4}(c).}}

In Section \ref{begin} we establish \whp\ that
$D_s'(t)=0$ for $t <\log^3n$ and $s \le r-2$.
We can use \eqref{Dst} to show that $\E D_s'(t)=o(1)$ for
$t\ll \t_{r-s}$ or
$t\geq t_{{s}}+\om n$. In which case the first moment method suffices.

We now deal with the vertices in $\bcN$.
Equation \eqref{ignore} shows that \whp\ $\bcN \seq \cB(t)$
for $t\geq t^*=10 \e_1 n \log n$.

and proves part (c).
\proofend
}

\subsection{Proof of Theorem \ref{th2}}
We combine Lemma \ref{lem2} and Theorem \ref{th4} with the results of
Molloy and Reed \cite{MR1,MR2}. We summarize
what we need from these two papers:
\begin{theorem}\label{MR}
Let $\l_0,\l_1,\ldots,\l_r\in [0,1]$ be such that $\l_0+\l_1+\cdots+\l_r=1$.
Suppose that $\bd=(d_1,d_2,\ldots,d_N)$
satisfies $|\set{j:d_j=s}|=(1+o(1))\l_sN$ for $s=0,1,\ldots,r$. Let $G_{n,\bd}$ be
chosen randomly from graphs with
vertex set $[N]$ and degree sequence $\bd$. Let
$$L=\sum_{s=1}^rs(s-2)\l_s.$$
\begin{description}\label{Leq0}
\item[(a)] If $L<0$ then \whp\ $G_{n,\bd}$ is sub-critical.
\item[(b)] If $L>0$ then \whp\ $G_{n,\bd}$ is super-critical.
Furthermore the unique giant component has
size $\th n$ where $\th$ is defined as follows: Let $\La=\sum_{s=1}^rs\l_s$.
Define $\a$ to be the smallest
positive solution to
\beq{La}
\La-2\a-\sum_{s=1}^rs\l_s\brac{1-\frac{2\a}{\La}}^{s/2}=0.
\eeq
Then
\beq{ff3}
\th=1-\sum_{s=0}^r\l_s\brac{1-\frac{2\a}{\La}}^{ s/2}.
\eeq
\end{description}
\end{theorem}
\vspace{-.3in}
\proofend

We now evaluate $L$ in the context of $\G(t)$.
Then Theorem \ref{th4}
implies that we can take
$$L=\sum_{s=0}^r\binom{r}{s}p_t^s(1-p_t)^{r-s}s(s-2)=rp_t((r-1)p_t-1).$$
Thus the critical value for $t$ is the one that gives $L=0$ and $p_t=\frac{1}{r-1}$.
One can easily check that this is
indeed the case
for $t^*$ as defined in \eqref{t*}.
{\mcyan Intuitively, this can be seen as follows. When we grow  components
outward from a given vertex, the  branching factor is $(r-1)p_t$, so when $(r-1)p_t<1$,
all components are of finite size.}

Parts (i) and (ii) of Theorem {\ref{th2}} follow. Here is a quick check on the claimed
size of $C_1(t)$. First of all $\La=rp_t>0$. We divide \eqref{La} by $\La$ and then
let $\f=1-\frac{2\a}{\La}$ so that \eqref{La} becomes \eqref{ff2}. With this value for $\f$
we see that \eqref{ff3} becomes \eqref{ff1} and then \whp\ $|C_1(t)|\sim \th N_t$ as claimed.

To prove (iii) we use the result of Hatami and Molloy \cite{HM} that when $|L|=O(n^{-1/3})$
the size of the giant is $n^{2/3+o(1)}$. At each step, {at most $r$ vertices in
$\G(t)$ change their degree and so
$L$ changes by $O(1/n)$. $L$ starts out at $r(r-2)$ and so at some time
it becomes equal to $O(n^{-1/3})$. This will happen at
$t\sim t^*$ and the conditions of \cite{HM} will be
satisfied. At this point \whp\ there are
$\Theta(n)$ vertices in $\G(t)$ and (iii) follows.
}
\proofend

\subsection{Proof of Theorem \ref{th3}}
In this section we study the number of components of a given size. In
principle one should be able to work
this out from Lemma \ref{lem2} and Theorem \ref{th4}. This has proven
more difficult than we anticipated.
Instead, we try to estimate the number directly. We can use these
lemmas though to argue that almost all
small components are trees. Indeed if we fix $t$ and condition on the values $D_s=D_s(t)$ satisfying
Theorem \ref{th4} then we have the following:
\begin{lemma}\label{not-trees}\

\vspace{-.25in}
\begin{description}
\item[(a)] If $t\leq \e_1n\log n$ then \whp\ there are at most
$n^{3\e_1}$ components of size $k\leq \e_1\log n$
that are not trees.
\item[(b)] If $t\geq \e_1n\log n$ then \whp\ there are no components
of size $k\leq \log^2n$ that are not trees.
\end{description}
\end{lemma}

\vspace{-.25in}
\proofstart
Let $N=|\cR(t)|$. Applying Lemma \ref{lem2} we see that the
expected number of sets of $k$ vertices that contain
at least $k$ edges is bounded by
$$\binom{N}{k}\binom{\binom{k}{2}}{k}\bfrac{r}{n}^k\leq \bfrac{rNe^2}{2n}^k.$$
To prove (a) we take $N=n$ and apply the Markov inequality. To prove (b) we take $N=N_t$.
\proofend

With this in mind we concentrate on the number of tree components
of size $k$ for some $k\leq \e_1\log n$. Since there
are \whp\ at most $n^{2\e_1}$ vertices that are not nice, we will
concentrate on counting the number of components
that are made up of nice vertices only. We will also assume that $t\leq t_0$, see \eqref{Dst}.

The following is Lemma 4 of  \cite{BFM}.
\begin{lemma}\cite{BFM}\label{smallbits0}
Let $\cT$ be  a  labeled  infinite $r$-regular tree.
Let $b_k$ be the number of labeled subtrees of size $k$
rooted at  vertex $v$ of $\cT$.
Then
\[
b_k = \frac{r}{(r-2)k+2}\binom{(r-1)k}{k-1}.
\]
\end{lemma}
\proofend

{Lemma  \ref{smallbits0}  counts labeled sub-trees. Let $T$ be such a
tree. Each edge $e=(x,y)$, with $x$ closest to $v$, is associated
with a label $\l_e$ indicating that it is the $\l_e$-th edge incident with $x$, in numerical order.}
Now consider the situation described in Lemma \ref{walkconfig}.
Fix {$v\in \cR(t)\cap\cN$}.
Assuming $v\in\cN$ introduces some conditioning on the allowable pairings
described in Lemma \ref{walkconfig}. However $\Pr(v\in \cN)=1-O(n^{2\e_1-1})$ (see \eqref{nice})
and then we can use
$$\frac{\Pr(A)-\Pr(\neg B)}{\Pr(B)}\leq \Pr(A\mid B)\leq \frac{\Pr(A)}{\Pr(B)}$$
for events $A,B$ to correct \eqref{ppt} below for this conditioning.
Consider the $k$ neighbourhood
of $v$ in the multi-graph on $[n]$
induced by a random pairing on $R(t)$. It is a tree. Now delete any edge that corresponds
to an edge $(x,y)$ with $x\in \cR(t),y\in\cB(t)$. Let $T$ be the
component that contains $v$. If $T$ has $k$ vertices then
$T$ corresponds to a tree component of $\G(t)$ with $k$ vertices.

So, fix a tree $T^*$ as counted in Lemma \ref{smallbits0} and let
us determine the probability that $T=T^*$.
{The probability space for this calculation is as follows: Let $F_t,R_t$ be
as in Lemma \ref{walkconfig}.
We have paired up the elements of $F_t$ and we are now considering random pairings of $R_t$.
{\mcyan To extend $F_t$ to $F$, we can generate the pairing of $R_t$ in any order we please.
Thus we start by pairing elements of $W_v$ the root of our tree.

 If we choose
$x\in W_v$ then the probability it is paired with $y\in W_z$, for some
$z\in \cR(t)$, is $\nu=\frac{r|\cR(t)|-1}{|R_t|-1}$.
Using Lemma \eqref{key}(a),(b), with $|\cR(t)|=\ooi N_t$ and $|R_t|=2|\cU(t)|$
we obtain
\beq{ppt}
\nu=\frac{r|\cR(t)|-1}{|R_t|-1}=\ooi e^{-(r-2)t/(\r r n)}=\ooi p_t,
\eeq
where $p_t$ is given by \eqref{pt}, and the $o(1)$ term is $o(\log^{-K}n)$
for any positive constant $K$. }

Suppose now that we have generated $O(\log n)$ pairings.
Then both $|\cR(t)|$ and $|R_t|$ change by $O(\log n)$ and, provided $t \le t_0(1-\e)$, $\e>0$,
since they are both of size $\Omega(n)$.
Choosing an unpaired
$x'\in W_{v'}$ for $v'\in \cR(t)$ we see that
the probability of $x$ being paired with some $y'\in W_{z'}$ where $z'\in \cR(t)$,
 is again $(1+o(1))p_t$.

To estimate $\Pr(T=T^*)$ we start at the root $v$ and examine the points paired with
$W_v$, where $W_v\subseteq R_t$, because $v\in \cR(t)$. The count in Lemma \ref{smallbits0}
assumes an ordering of the neighbours of each vertex
and by implication an ordering of $W_v$ and a statement about which members of $W_v$ are paired with
$W_{\cR(t)}$ and which should not. Suppose we pair $W_v$ with points
from $W_{x_i},i=1,2,\ldots,d$. Then we continue by
pairing up $W_{x_1}$ and then $W_{x_2}$ and so on. The factor $p_t^{k-1}$
is from the $k-1$ times we have to pair with
$W_{\cR(t)}$ and the factor $(1-p_t)^{(r-2)k+2}$ is from the number
of times we do not.

}
It follows that
\beq{TT*}
\Pr(T=T^*)=(1+\ole)p_t^{k-1}(1-p_t)^{(r-2)k+2}.
\eeq

It follows from \eqref{TT*} that
$$\E(N(k,t))=(1+\ole)N_t\frac{b_k}{k}p_t^{k-1}(1-p_t)^{(r-2)k+2}.$$
It remains to prove concentration around the mean. We use the Chebyshev
inequality. We fix two vertex disjoint
trees $T_1,T_2$ in $G$. Arguing as above we see that
$$\Pr(T_1,T_2\text{ are components of }\G(t))\leq (1+o(1))\prod_{i=1}^2
\Pr(T_i\text{ is a component of }\G(t)).$$
Provided $\E(N(k,t))\to\infty$ and it does so for
$(\t_{k(r-2)+2})^{1+\d} \le t \le (1-\e)t_{k}$ we can use the Chebychev
inequality to prove Theorem \ref{th3}(a) for any fixed $t$. Note that here we have used concentration
in the configuration model to imply concentration in the simple graph model.

{\grey
We now prove Theorem \ref{th3}(b).
For $k$ constant, $\E(N(k,t))=(1+o(1))\eta(k,t)\ge n^{\e'}$,  throughout the range $\e n \le t \le (1-\e)t_{k-1}$.
Concentration can be established by the methods used in Section \ref{proof of th3}.
Let $A$ be a large constant. For $\e n \le t \le An$ we use the
Chebychev inequality directly, and for $An \le t \le (1-\e)t_{k-1}$
we use the interpolation method. For  $\e n \le t \le An$, $N(k,t) \ge a_k n$
where $a_k$ is some constant. Using \eqref{toterr}, for some $\d>0$,
\[
\Pr(|N(k,t)-\E N(k,t)| \ge \E N(k,t) o(n^{-\d})) =o(n^{-\d}).
\]
Interpolate $\e n \le t \le An$ at $h=n^{\d/2}$ points $(t_1,...,t_h)$, $\ell=(A-\e)n^{1-\d/2}$
apart. Let $t \in (t_j,t_{j+1})$. Then $|N(k,t)-N(k,t_j)| \le r\ell$ and 
$|\eta(k,t)-\eta(k,t_j)|=O(\ell)$. It follows that, with probability $1-o(n^{-\d/2})$,
$|N(k,t)-\E N(k,t)| \le \E N(k,t) o(n^{-\d/2})$ for all $\e n \le t \le An$.
For $AN \le t \le (1-\e)t_{k-1}$ we use the interpolation method of Section \ref{proof of th3} as follows.
\Whp\ the maximum component size is $O(\log n)$ and there are $O(\log n)$ vertices on non-tree
components. Let $M(k,t)$ be the number of vertices on components of
size at least $k$. Then $M(k,t)$ is monotone non-increasing and
\begin{eqnarray*}
M(k,t)&=& \sum_{j = k}^{O(\log n)} jN(j,t)+O(\log n)\\
k N(k,t)&=& M(k-1,t)-M(k,t)+O(\log n).
\end{eqnarray*}
Applying the interpolation method to $M(k,t)$ will complete the proof of
Theorem \ref{th3}(b).
}

\subsection{In the beginning}\label{begin}
{We remind the reader that our proofs so far, we have assumed $t\ge \log^3n$.
We now consider the first few moves of the walk.}
Using Lemma \ref{lem2} we see that for $1\leq t\leq \log^3n$ we have
that $\G(t)$ is a random graph with a degree sequence
of the following form: There are $n-s$ vertices of degree $r$, where
$s\leq rt$, and $s$ vertices of degree $<r$.
If the minimum degree in $\G(t)$ is at least one then \whp\ we find
that $\G(t)$ is connected. Indeed, let $V_r$ be
the set of vertices of degree $r$ in $\G(t)$. We argue that \whp
\begin{align}
&V_r\text{ induces a connected subgraph of }\G(t).\label{a}\\
&\text{Each $x\in \cR(t)\setminus V_r$ is adjacent to $V_r$}.\label{b}
\end{align}
For $k$ even let
$$\f(k)=\frac{k!}{(k/2)!2^{k/2}}$$
be the number of ways of partitioning $[k]$ into $k/2$ pairs.

Let $m=O(\log n)$ be the sum of the degrees, in $\G(t)$, of the
vertices in $\cR(t)\setminus V_r$. Then
\begin{align}
\Pr(\eqref{a} fails)&\leq \sum_{k=3}^{n/2}\sum_{l=0}^m\binom{n}{k}\binom{m}{l}
\frac{\f(kr+l)\f((n-s-k)r+m-l)}
{\f((n-s)r+m)}\label{afail}\\
&=\sum_{k=3}^{n/2}\sum_{l=0}^m\binom{n}{k}\binom{m}{l}\frac{\binom{(r(n-s)+m)/2}{(kr+l)/2}}{\binom{r(n-s)+m}{kr+l}}\nonumber\\
&\le\sum_{k=3}^{n/2}\sum_{l=0}^m\binom{n}{k}\binom{m}{l}\frac{1}{\binom{(r(n-s)+m)/2}{(kr+l)/2}}\nonumber\\
&\leq \sum_{k=3}^{n/2}\sum_{l=0}^m\binom{n}{k}\binom{m}{l}\bfrac{kr+l}{r(n-s)+m}^{(kr+l)/2}
\nonumber\\
&=o(1).\nonumber
\end{align}
{\bf Explanation of \eqref{afail}}: Choose a set of $k$ vertices
$S$ of degree $r$ and $l$ points from the
$m$ points in $W$ associated with vertices $T$
of degree less than $r$. (Some of the points associated with these
latter vertices have already been paired). Now pair up the $kr+l$
points randomly and the remaining points
randomly. The $l$ points contain the edges between $S$ and $T$.

The probability that \eqref{b} fails is $O(\log n/n)$.
A vertex of $\cR(t)\setminus V_r$ of degree $d$
has an $O(n^{-d})$
chance of not being connected to $V_r$.

So we only have to deal with the possibility that there
are isolated vertices in $\G(t)$ for $t\leq \log^3n$.
So consider the event
$$\cA(t)=\set{\exists v\in N(W_t):\;v\in \cR(t)\text{ and }N(v)\subseteq \cB(t)}.$$
We claim that
\beq{cA}
\Pr(\cA(t))=O\bfrac{\log^3n}{(r-1)^{\ell_1/2}}.
\eeq
It follows that
$$\Pr(\eqref{b}\ fails)\leq \sum_{t=1}^{\log^3n}\Pr(\cA(t))=o(1).$$
To prove \eqref{cA} fix $t$ and a neighbour $v$ of $W_t$. Equation \eqref{nice1} implies that
there is at least one neighbour
$w$ of $v$ that is not contained in a small cycle. If
$w\neq X_t$ then to reach $w$ the walk $\cW$ must emulate
a walk on the infinite tree $\cT$ that starts at distance
$\ell_1/2$ from the root and visits it. This has
probability $1/(r-1)^{\ell_1/2}$ and this must be inflated
 by $\log^3n$ to account for $\log^3n$ possible starting times.
If $w=X_t$ then to visit another neighbour of $v$ then we must
first reach distance at least $\ell_1/2$
and then we can repeat the argument and use inequality \eqref{Hoef}.

\section{The evolution of $\G(t)$ in $G_{n,p}$ and $D_{n,p}$}\label{Gnp}

We first consider $G_{n,p}$ and prove Theorem \ref{th1}.

We first note some properties of the degree sequence $d_G(v)$ of $G_{n,p}$.
We assume that $c=n^{1/\om}$ where $\om\to\infty$.
Let $\om_1=\log^{1/3}n$.
For a fixed vertex $v$, its degree $d_G(v)$ satisfies
$$\Pr(|d_G(v)-c\log n|\geq \om_1(c\log n)^{1/2})\leq 2e^{-\om_1^2/3}.$$
This follows from Chernoff bounds on the tail of the binomial distribution.
Denote  by $\cN_d$, vertices of $G_{n,p}$ which have degrees in the range
$c\log n\pm \om_1(c\log n)^{1/2}$.
By the Markov inequality, we see that \whp\ $|\neg(\cN_d)|=(ne^{-\om_1^2/4})$.

Let $|\cR(t)|=N$.
Because $\G(t)$ has the distribution $G_{N,p}$, we only need
good estimates of $N$. We can get these using
Lemma \ref{First}.
Fix  a vertex $v \in \cN_d$.
It is shown in \cite{CFGiant} that
\whp\ $R_v=1+O(1/\log n)$ for all $v\in V$.
Let $\e=1/(\log n)^{1/6})$, then for $v \in \cN_d$, $\pi_v =(1+O(\e)) (1/n)$.
Thus $p_v=(1+O(\e))(1/n)$.

For $t \ge \log^3 n$,
 \eqref{atv1} implies that
$$
\Pr(v\in \cR(t))=(1+O(\log n/n))e^{-(1+O(\e))t/n}. 
$$
Recall that $t_\th=n(\log \log n+(1+\th) \log c)$, and
 assume that $t=t_\th$, where $\th=O(1)$ then
$$\E(|\cN_d\cap \cR(t)|)=(1+o(1))\frac{n}{c^{1+\th}\log n}.$$

Regarding concentration,
we can argue as in the proof \eqref{var} that if $v,w\in \cN_d$ are at
distance at least {$\om_1/2$} in $G$ then
\beq{fire}
\Pr(\cE_v\cap \cE_w)=(1+o(1))\Pr(\cE_v)\Pr(\cE_w).
\eeq
where $\cE_v=\set{v\in \cR(t)}$.

Then if $X=|\cR(t)\cap \cN_d|$ then
{\red
\begin{align*}
\E(X(X-1))&\leq \sum_{v\in \cN_d}\sum_{\substack{w\in\cN_d\\dist(v,w)\geq \om_1/2}}
\Pr(\cE_v\cap \cE_w)+
\sum_{v\in \cN_d}\sum_{\substack{w\in\cN_d\\dist(v,w)\leq \om_1/2}}\Pr(\cE_v\cap \cE_w)\\
&\leq (1+o(1))\sum_{v\in \cN_d}\sum_{\substack{w\in\cN_d\\dist(v,w)
\geq \om_1/2}}\Pr(\cE_v)\Pr(\cE_w)+
\sum_{v\in \cN_d}\sum_{\substack{w\in\cN_d\\dist(v,w)\leq \om_1/2}}\Pr(\cE_v)\\
&\leq (1+o(1))\E(X)^2+O((c\log n)^{\om_1/2})\E(X)
\end{align*}}
which implies that $\Var(X)=o(\E(X)^2)$ and then the Chebyshev inequality implies that
$X\sim \E(X)\sim \frac{n}{c^{1+\th}\log n}$ \whp.

The vertices outside $\cN_d$ only contribute $o(N)$ \whp\ and
thus
\[
N(t_\th)=\ooi \E N(t_\th) = \ooi n/(c^{1+\th} \log n).
\]
The threshold for the giant component in $G_{N,p}$ is at $Np=1$, i.e. as $\th \ra 0$ from below.
Theorem \ref{th1} follows immediately from this.

We next consider $D_{n,p}$ and prove Theorem \ref{th1dir}. As  the details are similar to those above, our discussion
will be brief.  If $np=c \log n$ and $(c-1) \log n \rai$ then \whp \ $D_{n,p}$ is strongly connected, so
a random walk on $D_{n,p}$ is ergodic.
It was established in \cite{CFDir} that if $(c-1) \log n \rai$, then
almost all vertices $v$ have
have stationary distribution $\pi_v=\ooi/n$. By the method of deferred decisions,
$\vec \G(t)$ is a random digraph $D_{N(t),p}$ on $N(t)=|\cR(t)|$ vertices,
where $N(t)$ is given as above. It was proved in \cite{Karp} that the threshold for the emergence of a giant
strongly connected component in $D_{n,p}$ is at $np \sim 1$. Theorem \ref{th1dir} follows.

\end{document}